\documentclass[a4paper, 11pt]{article}
\usepackage{amsmath, amssymb, amsthm}
\usepackage[margin=21mm]{geometry}
\usepackage{microtype}
\usepackage[colorlinks=true, linkcolor=blue, citecolor=blue, urlcolor=blue]{hyperref}

\newtheorem{theorem}{Theorem}[section]
\newtheorem{lemma}[theorem]{Lemma}

\theoremstyle{remark}
\newtheorem{remark}[theorem]{Remark}

\newcommand{\E}{\mathbb{E}}
\newcommand{\Prob}{\mathbb{P}}
\newcommand{\Hyp}{\mathrm{Hyp}}
\newcommand{\cS}{\mathcal S}
\DeclareMathOperator{\tr}{tr}
\DeclareMathOperator{\Var}{Var}

\title{\vspace{-1.2em}The number of $k$-factors of a random regular graph}
\author{%
Katerina Adler \qquad Simi Haber\\
{\normalsize Department of Mathematics, Bar-Ilan University, Ramat Gan, Israel}%
}
\date{\vspace{-0.8em}October 2026\vspace{-0.6em}}

\begin{document}
\maketitle

\begin{abstract}
Let $Z_k(G)$ be the number of $k$-factors ($k$-regular spanning subgraphs) of a uniformly random simple $d$-regular graph $G\sim G(n,d)$, where $d\ge 3$ and $1\le k\le d-1$ are fixed integers such that $dn$ and $kn$ are even. We determine the asymptotic second moment of $Z_k(G)$ and show that $Z_k(G)/\E Z_k(G)$ converges in distribution to an explicit random variable $W$ that depends on the short-cycle counts of $G$ but not on $k$. Moreover, for any $1\le k,\ell\le d-1$, the normalised counts on the same graph satisfy
\[
\frac{Z_k(G)}{\E Z_k(G)}-\frac{Z_\ell(G)}{\E Z_\ell(G)}\;\xrightarrow{\;p\;}\;0,
\]
and, since $W>0$ almost surely, their ratio converges to $1$ in probability. This extends results of Bollob\'as--McKay and Janson for $k=1$ and Robalewska for $k=2$ to all $1\le k\le d-1$. In the pairing model, the second moment reduces to a one-dimensional lattice sum over the overlap of two $k$-factors. We evaluate this sum using the Poisson-binomial representation of the hypergeometric distribution, exponential tilting, and a uniform local limit theorem.
\end{abstract}


\section{Introduction}\label{sec:intro}

Fix integers $d\ge 3$ and $1\le k\le d-1$, and let $n$ be a positive integer such that $dn$ and $kn$ are even. A \emph{$k$-factor} of a graph $G$ is a spanning $k$-regular subgraph. Let $G(n,d)$ be the probability space of uniformly random labelled simple $d$-regular graphs on $[n]=\{1,\dots,n\}$. Write $Z_k(G)$ for the number of $k$-factors of $G$, and $X_i(G)$ for the number of cycles of length $i$ in $G$. Throughout,
\begin{equation}\label{eq:lambda-delta}
\lambda_i=\frac{(d-1)^i}{2i},\qquad
\delta_i=\frac{(-1)^i}{(d-1)^i}\qquad(i\ge 1),
\end{equation}
so $1+\delta_i>0$ for all $i$ and $\sum_{i\ge1}\lambda_i\delta_i^2=\tfrac12\sum_{i\ge1}(d-1)^{-i}/i=\tfrac12\log\frac{d-1}{d-2}<\infty$. We prove that the short cycles of $G$ govern the relative fluctuations of $Z_k(G)$ around its mean, independently of $k$.

\begin{theorem}\label{thm:main}
Fix integers $d\ge 3$ and $1\le k\le d-1$. Let $n\to\infty$ along integers for which $dn$ and $kn$ are even, and let $G\sim G(n,d)$. Let $(\xi_i)_{i\ge 3}$ be independent $\mathrm{Poisson}(\lambda_i)$ random variables.
\begin{enumerate}
\item[\textup{(i)}] \textup{(mean)}
\begin{equation}\label{eq:mean-simple}
\E_{G(n,d)}[Z_k(G)]=\sqrt{2}\,e^{1/4}\Bigl(\binom{d}{k}\,k^{k/2}(d-k)^{(d-k)/2}d^{-d/2}\Bigr)^{n}\bigl(1+o(1)\bigr).
\end{equation}
\item[\textup{(ii)}] \textup{(second moment)}
\begin{equation}\label{eq:second-moment-simple}
\frac{\E_{G(n,d)}\bigl[Z_k(G)^2\bigr]}{\bigl(\E_{G(n,d)}[Z_k(G)]\bigr)^2}
\;\longrightarrow\;
\exp\Bigl(\sum_{i\ge 3}\lambda_i\delta_i^2\Bigr)
\;=\;
e^{-\frac{2d-1}{4(d-1)^2}}\sqrt{\frac{d-1}{d-2}} .
\end{equation}
\item[\textup{(iii)}] \textup{(limit distribution)} The infinite product
\begin{equation}\label{eq:W-def}
W\;:=\;\prod_{i=3}^{\infty}(1+\delta_i)^{\xi_i}\,e^{-\lambda_i\delta_i}
\end{equation}
converges almost surely and in $L^2$. It satisfies $W>0$ almost surely, $\E W=1$, and $\E W^2=\exp\bigl(\sum_{i\ge3}\lambda_i\delta_i^2\bigr)$. Moreover,
\begin{equation}\label{eq:limit-distribution}
\frac{Z_k(G)}{\E_{G(n,d)}[Z_k(G)]}\;\Longrightarrow\;W,
\end{equation}
jointly with $(X_3(G),\dots,X_m(G))\Rightarrow(\xi_3,\dots,\xi_m)$ for every fixed $m$.
\item[\textup{(iv)}] \textup{(linkage)} For any $1\le k,\ell\le d-1$, let $n$ tend to infinity along integers for which $dn$, $kn$, and $\ell n$ are even. The normalised counts on the same graph $G\sim G(n,d)$ satisfy
\begin{equation}\label{eq:linked}
\frac{Z_k(G)}{\E_{G(n,d)}[Z_k(G)]}-\frac{Z_\ell(G)}{\E_{G(n,d)}[Z_\ell(G)]}
\;\xrightarrow{\;p\;}\;0,
\end{equation}
and consequently $\dfrac{Z_k(G)/\E_{G(n,d)}[Z_k(G)]}{Z_\ell(G)/\E_{G(n,d)}[Z_\ell(G)]}\xrightarrow{\;p\;}1$. Here the ratio is set equal to $1$ on the event $\{Z_\ell(G)=0\}$, whose probability tends to $0$.
\end{enumerate}
\end{theorem}

\paragraph{Previous work.}
For $k=1$ (perfect matchings), Bollob\'as and McKay~\cite{bollobas_mckay1986} determined the asymptotics of the first two moments of $Z_1(G)$. Janson~\cite{janson1995} obtained the limit distribution of $Z_1(G)/\E Z_1(G)$ by the small subgraph conditioning method. Robinson and Wormald introduced this method for Hamilton cycles in cubic graphs~\cite{robinson_wormald1992} and applied it to perfect matchings~\cite{robinson_wormald1994}. For $k=2$, Robalewska~\cite{robalewska1996} computed the second moment and established the limit distribution.

The complement of a $k$-factor in a $d$-regular graph is a $(d-k)$-factor. Thus $Z_k(G)=Z_{d-k}(G)$, and the cases $k=d-1$ and $k=d-2$ follow. Theorem~\ref{thm:main} was therefore known for $k\in\{1,2,d-2,d-1\}$. In particular, the normalised numbers of perfect matchings and $2$-factors have the same limit distribution~\cite[Corollary~4]{robalewska1996}. Wormald also observed that they are asymptotically linked: $Z_1(G)/\E Z_1(G)\sim Z_2(G)/\E Z_2(G)$ asymptotically almost surely~\cite[Corollary~4.10 and the remark following it]{wormald1999}. Part~(iv) extends this result to all pairs $k,\ell$.

The remaining cases $3\le k\le d-3$ occur when $d\ge 6$. For these $k$, the contiguity of $G(n,d)$ and the edge-disjoint superposition $G(n,k)\oplus G(n,d-k)$ follows from~\cite[Corollary~4.17]{wormald1999}, which combines the contiguity results collected in~\cite[Theorem~4.15]{wormald1999}; see also~\cite{janson1995}. The first moment in the pairing model is elementary; see~\eqref{eq:EZk-pairing} below. For fixed $d$, we are not aware of a previous determination of the second moment or the limit distribution of $Z_k(G)$ for $3\le k\le d-3$. Isaev et al.~\cite[Eq.~(4.2)]{isaev2025} prove the bound $\E[Z_k(G)^2]\le(1+o(1))(\E Z_k(G))^2$ when $k^2\bigl(n/\min\{d-k,\,n-d+k\}\bigr)^{2k}\le\tfrac13\log n$. In the growing-degree regime, Gao~\cite{gao2023} proved the asymptotic normality of $Z_1(G)$ for $d\to\infty$ with $d=O(n^{1/7}/\log n)$. Gao conjectured this result for all $Z_k(G)$ when $\min(d,n-d)\to\infty$.

For fixed $d$ and $k$, the mean formula in part~(i), including the constant $\sqrt{2}\,e^{1/4}$, follows from the enumeration of Hasheminezhad and McKay~\cite{hasheminezhad_mckay2023}. It equals the number of partitions of the complete graph into spanning regular subgraphs of degrees $n-1-d$, $k$ and $d-k$, divided by the number of $d$-regular graphs. The second moment and the limit distribution are proved here.

\paragraph{Contribution.}
We prove Theorem~\ref{thm:main} by Janson's small subgraph conditioning method~\cite{janson1995} in the pairing model. The proof uses two lemmas. Lemma~\ref{lem:short-cycles} proves the joint factorial moments of $Z_k$ and the short-cycle counts. A $2\times2$ transfer matrix with spectrum $\{d-1,-1\}$, independent of $k$, explains the universality in $k$. Janson states the same matrix and the resulting Poisson means without proof in his section on contiguity and asymptotic distributions~\cite{janson1995}.

Next, Lemma~\ref{lem:second-moment} proves the second moment of $Z_k$. The overlap of two $k$-factors is a sum of i.i.d.\ hypergeometric random variables, one for each vertex, so the second moment reduces to a one-dimensional lattice sum over this overlap (Lemma~\ref{lem:reduction}). We evaluate this sum by a Laplace-type argument. The key step is that a hypergeometric distribution is a sum of independent Bernoulli variables with distinct parameters (Lemma~\ref{lem:poisson-binomial}; equivalently, its generating polynomial has only real zeros~\cite{vatutin_mikhailov1982,pitman1997}). This makes the tilted family explicit, gives a strictly concave rate-function comparison (Lemma~\ref{lem:rate}), and yields a uniform local limit theorem for the tilted sums (Lemma~\ref{lem:llt}). The passage from the pairing model to simple graphs is standard~\cite[Lemma~4.4]{wormald1999}. For completeness, Section~\ref{sec:proof} derives the conditional first and second moments from uniform integrability.

\paragraph{Organisation.}
Section~\ref{sec:pairing} recalls the pairing model, the first moment, and Janson's small subgraph conditioning theorem. Section~\ref{sec:short-cycles} treats short cycles, Section~\ref{sec:second-moment} treats the second moment, and Section~\ref{sec:proof} proves Theorem~\ref{thm:main}.

\section{The pairing model and small subgraph conditioning}\label{sec:pairing}

We use the pairing model $\mathcal P_{n,d}$~\cite{bollobas1980,wormald1999}. Each vertex $v\in[n]$ is represented by a cell $V_v$ of $d$ labelled points, and a pairing $P\in\mathcal P_{n,d}$ is a perfect matching of the $dn$ points. The number of perfect matchings on $2\ell$ points is
\[
M(2\ell)=\frac{(2\ell)!}{\ell!\,2^{\ell}}=(2\ell-1)!!\qquad(\ell\ge0),
\]
where $M(0)=1$; set $M(x)=0$ for odd $x$.

Let $\Prob$ and $\E$ denote probability and expectation under the uniform distribution on $\mathcal P_{n,d}$. Contracting each cell $V_v$ to a vertex projects $P$ to a $d$-regular multigraph $G(P)$. A pair within a cell is a loop, and two pairs joining the same two cells are parallel edges. Let $X_i=X_i(P)$ be the number of $i$-cycles of $G(P)$. Here $X_1$ is the number of loops, $X_2$ is the number of unordered pairs of parallel edges, and, for $i\ge3$, $X_i$ is the number of cycles of length $i$. With these conventions, for every fixed $m$,
\begin{equation}\label{eq:poisson-cycles}
(X_1,\dots,X_m)\;\Longrightarrow\;(\xi_1,\dots,\xi_m),\qquad \xi_i\sim\mathrm{Poisson}(\lambda_i)\ \text{independent}
\end{equation}
\cite[Theorem~2.6]{wormald1999}; see also the proof of Theorem~2 in~\cite{bollobas1980}. The multigraph $G(P)$ is simple precisely on the event $\cS:=\{X_1=X_2=0\}$. The probability of this event tends to $e^{-\lambda_1-\lambda_2}=e^{-(d^2-1)/4}>0$. Conditional on $\cS$, the graph $G(P)$ is uniformly distributed over simple $d$-regular graphs; that is, $G(P)\sim G(n,d)$.

A \emph{$k$-factor of $P$} is a subset $F\subseteq P$ of pairs that covers exactly $k$ points in each cell $V_v$. Let $Z_k(P)$ denote the number of $k$-factors of $P$. If $G(P)$ is simple, the pairs of $P$ correspond bijectively to the edges of $G(P)$, and the $k$-factors of $P$ correspond to those of $G(P)$. Thus $Z_k(P)=Z_k(G(P))$ on $\cS$.

A $k$-factor $F$ determines the sets $A_v\subseteq V_v$ of points that it covers, where $|A_v|=k$. There are $\binom{d}{k}^n$ choices for these sets, $M(kn)$ ways to pair the $kn$ chosen points to form $F$, and $M((d-k)n)$ ways to pair the remaining $(d-k)n$ points to form $R=P\setminus F$. Hence the exact first moment is
\begin{equation}\label{eq:EZk-pairing}
\E Z_k(P)
\;=\;
\binom{d}{k}^n\,\frac{M(kn)\,M\bigl((d-k)n\bigr)}{M(dn)}
\;=\;
\sqrt{2}\,\Bigl(\binom{d}{k}\,k^{k/2}(d-k)^{(d-k)/2}d^{-d/2}\Bigr)^{n}\bigl(1+O(n^{-1})\bigr),
\end{equation}
where the asymptotic form follows from Stirling's formula in the form $M(2\ell)=\sqrt2\,(2\ell/e)^{\ell}(1+O(\ell^{-1}))$. The powers of $n$ and of $e$ cancel because $kn+(d-k)n=dn$.

The same count describes the uniform distribution on pairs $(P,F)$ in which $F$ is a $k$-factor of $P$. The sets $(A_v)_{v\in[n]}$ are independent uniform $k$-subsets of the cells. Conditional on these sets, $F$ is a uniform perfect matching of the chosen points and $R$ is an independent uniform perfect matching of the remaining points. Let $\E_\Omega$ denote expectation for a uniformly random pair $(P,F)$ in $\Omega:=\{(P,F):F\text{ is a $k$-factor of }P\}$. For every function $\Psi$ of $P$,
\begin{equation}\label{eq:weighted}
\frac{\E[Z_k(P)\,\Psi(P)]}{\E Z_k(P)}\;=\;\E_{\Omega}\bigl[\Psi(P)\bigr].
\end{equation}

We use the following form of the small subgraph conditioning theorem. It is a special case of Theorem~1 of~\cite{janson1995}, with his second condition replaced by the factorial-moment condition in the lemma following that theorem. It also follows from Theorems~4.1 and~4.3 of~\cite{wormald1999}. Here $(x)_r=x(x-1)\cdots(x-r+1)$ denotes the falling factorial.

\begin{theorem}[small subgraph conditioning]\label{thm:janson}
Let $\lambda_i>0$ and $\delta_i>-1$ $(i\ge1)$ be constants. For each $n$, let $Y_n\ge0$ and $X_{1n},X_{2n},\dots$ be random variables on a common probability space, where the $X_{in}$ are nonnegative and integer valued and $0<\E Y_n<\infty$. Suppose that
\begin{enumerate}
\item[\textup{(i)}] for each fixed $m$, $(X_{1n},\dots,X_{mn})\Rightarrow(\xi_1,\dots,\xi_m)$, where the $\xi_i\sim\mathrm{Poisson}(\lambda_i)$ are independent;
\item[\textup{(ii)}] for each fixed $m$ and all nonnegative integers $r_1,\dots,r_m$,
$\displaystyle \frac{\E\bigl[Y_n\prod_{i=1}^m(X_{in})_{r_i}\bigr]}{\E Y_n}\longrightarrow\prod_{i=1}^m\bigl(\lambda_i(1+\delta_i)\bigr)^{r_i}$;
\item[\textup{(iii)}] $\sum_{i\ge1}\lambda_i\delta_i^2<\infty$;
\item[\textup{(iv)}] $\E Y_n^2/(\E Y_n)^2\longrightarrow\exp\bigl(\sum_{i\ge1}\lambda_i\delta_i^2\bigr)$.
\end{enumerate}
Then the product $W^{\mathrm{pair}}:=\prod_{i\ge1}(1+\delta_i)^{\xi_i}e^{-\lambda_i\delta_i}$ converges almost surely, and, for each fixed $m$,
\[
\Bigl(\frac{Y_n}{\E Y_n},\,X_{1n},\dots,X_{mn}\Bigr)\;\Longrightarrow\;\bigl(W^{\mathrm{pair}},\,\xi_1,\dots,\xi_m\bigr)
\qquad\text{jointly.}
\]
\end{theorem}

We need the following elementary properties of the limiting products for the constants $\lambda_i,\delta_i$ in~\eqref{eq:lambda-delta}.

\begin{lemma}\label{lem:W-properties}
Let $\lambda_i,\delta_i$ be as in~\eqref{eq:lambda-delta}, let $(\xi_i)_{i\ge1}$ be independent $\mathrm{Poisson}(\lambda_i)$ variables, and let $I\subseteq\{1,2,\dots\}$. Then the partial products $W_I^{(m)}:=\prod_{i\in I,\,i\le m}(1+\delta_i)^{\xi_i}e^{-\lambda_i\delta_i}$ converge as $m\to\infty$, almost surely and in $L^2$, to a limit $W_I$ such that $W_I>0$ almost surely, $\E W_I=1$, and $\E W_I^2=\exp\bigl(\sum_{i\in I}\lambda_i\delta_i^2\bigr)$. If $I\cap J=\emptyset$, then $W_I$ and $W_J$ are independent.
\end{lemma}

\begin{proof}
Since $|\delta_i|\le(d-1)^{-i}\le\tfrac12$, we have $|\log(1+\delta_i)-\delta_i|\le\delta_i^2$ and $|\log(1+\delta_i)|\le2|\delta_i|$. The summands $\xi_i\log(1+\delta_i)-\lambda_i\delta_i$ of $\log W_I^{(m)}$ are independent. Their means are $\lambda_i(\log(1+\delta_i)-\delta_i)=O(\lambda_i\delta_i^2)$, and their variances are $\lambda_i\log^2(1+\delta_i)=O(\lambda_i\delta_i^2)$. Both are summable because $\sum_{i\ge1}\lambda_i\delta_i^2<\infty$, as noted after~\eqref{eq:lambda-delta}. Kolmogorov's two-series theorem therefore shows that the series converges almost surely and that $W_I=\exp(\lim_m\log W_I^{(m)})\in(0,\infty)$ almost surely.

Each factor $(1+\delta_i)^{\xi_i}e^{-\lambda_i\delta_i}$ has mean $e^{\lambda_i\delta_i}e^{-\lambda_i\delta_i}=1$ and second moment $e^{\lambda_i((1+\delta_i)^2-1)}e^{-2\lambda_i\delta_i}=e^{\lambda_i\delta_i^2}$. Thus $(W_I^{(m)})_m$ is a nonnegative martingale with bounded second moments, $\E(W_I^{(m)})^2=\exp(\sum_{i\in I,i\le m}\lambda_i\delta_i^2)$. It therefore converges in $L^2$, and $\E W_I=1$, $\E W_I^2=\exp(\sum_{i\in I}\lambda_i\delta_i^2)$. Independence for disjoint index sets follows from the independence of the $\xi_i$.
\end{proof}

\section{Short cycles}\label{sec:short-cycles}

\begin{lemma}[short-cycle moments]\label{lem:short-cycles}
Fix $d\ge 3$ and $1\le k\le d-1$. For any fixed $m\ge1$ and nonnegative integers $r_1,\dots,r_m$,
\begin{equation}\label{eq:factorial-moments}
\frac{\E\Bigl[Z_k(P)\prod_{i=1}^m (X_{i})_{r_i}\Bigr]}{\E Z_k(P)}
\;\longrightarrow\;
\prod_{i=1}^m \bigl(\lambda_{i}(1+\delta_{i})\bigr)^{r_i}
\;=\;\prod_{i=1}^m\Bigl(\frac{(d-1)^i+(-1)^i}{2i}\Bigr)^{r_i}.
\end{equation}
\end{lemma}

\begin{proof}
By~\eqref{eq:weighted}, the left-hand side equals $\E_\Omega\prod_i(X_i)_{r_i}$ for a uniform pair $(P,F)\in\Omega$. Generate $(P,F)$ as described after~\eqref{eq:EZk-pairing}: choose independent uniform $k$-subsets $A_v\subseteq V_v$, then independent uniform perfect matchings $F$ of $\mathcal A:=\bigcup_vA_v$ and $R$ of its complement, and set $P=F\cup R$. Call the points in $\mathcal A$ \emph{$F$-points} and the remaining points \emph{$R$-points}. Every pair of $P$ joins either two $F$-points (an \emph{$F$-pair}) or two $R$-points (an \emph{$R$-pair}).

\emph{Directed closed walks.}
By dihedral symmetry, for every $i\ge1$ the number of directed closed walks of length $i$ in $P$ is $2i X_i$. Such a walk consists of a sequence $(v_1,\dots,v_i)$ of distinct vertices together with points $y_j,x_j\in V_{v_j}$, $y_j\ne x_j$, such that $\{x_j,y_{j+1}\}\in P$ for each $j$, indices taken modulo $i$. The walk arrives at $v_j$ through $y_j$ and leaves through $x_j$.

\emph{A single cycle.}
Fix distinct vertices $v_1,\dots,v_i$ and a type $\tau_j\in\{F,R\}$ for each edge $j$, that is, for the pair $\{x_j,y_{j+1}\}$. If $\tau_j=F$, the points $x_j$ and $y_{j+1}$ must both be $F$-points; if $\tau_j=R$, they must both be $R$-points. At vertex $v_j$, the entering edge $j-1$ and the leaving edge $j$ have types $(\tau_{j-1},\tau_j)$. The number of ordered pairs $(y_j,x_j)$ of distinct points of $V_{v_j}$ with the prescribed memberships is $w(\tau_{j-1},\tau_j)$, where
\[
w(F,F)=k(k-1),\qquad w(R,R)=(d-k)(d-k-1),\qquad w(F,R)=w(R,F)=k(d-k).
\]
Conditional on the sets $A_v$, the $i$ prescribed pairs are pairwise disjoint. If $a$ of them have type $F$ and $b=i-a$ have type $R$, the independence of $F$ and $R$ gives the probability that all of them belong to $P$.
\begin{equation}\label{eq:pair-prob}
\begin{aligned}
\frac{M(kn-2a)}{M(kn)}\cdot\frac{M((d-k)n-2b)}{M((d-k)n)}
&=\prod_{j=0}^{a-1}\frac{1}{kn-1-2j}\;\prod_{j=0}^{b-1}\frac{1}{(d-k)n-1-2j}\\
&=(kn)^{-a}\bigl((d-k)n\bigr)^{-b}\bigl(1+O(n^{-1})\bigr).
\end{aligned}
\end{equation}
Each $F$-edge has pairing probability $\sim 1/(kn)$ among the $kn$ points assigned to $F$, and each $R$-edge has pairing probability $\sim 1/((d-k)n)$ among the remaining points. There are $(n)_i=n^i(1+O(n^{-1}))$ choices for the vertices. At a vertex with entering type $\tau$ and leaving type $\tau'$, the weight $w(\tau,\tau')$ and the pairing probability of the leaving edge contribute the entry $T_{\tau\tau'}$, equal to $w(\tau,\tau')/k$ or $w(\tau,\tau')/(d-k)$ according as $\tau'$ is $F$ or $R$. Recording the sequence of edge types around the cycle by this transition matrix, and summing over the $2^i$ type sequences, shows that the expected number of directed closed walks of length $i$ is
\[
\E_\Omega=(1+O(n^{-1}))\tr(T^i).
\]
Order the types as $(F,R)$, so the first row and column refer to $F$-edges and the second to $R$-edges. Then
\[
T=\begin{pmatrix} k-1 & k \\ d-k & d-k-1 \end{pmatrix},\qquad
\tr T=d-2,\qquad \det T=-(d-1),
\]
so the eigenvalues of $T$ are $d-1$ and $-1$ for every $k$. (With $D=\mathrm{diag}(\sqrt{k},\sqrt{d-k})$, the matrix $D^{-1}TD$ is symmetric, with diagonal entries $k-1$, $d-k-1$ and off-diagonal entries $\sqrt{k(d-k)}$.) Hence $\tr(T^i)=(d-1)^i+(-1)^i$. Since the number of these walks is $2i X_i$,
\[
\E_\Omega X_i\;\longrightarrow\;\frac{(d-1)^i+(-1)^i}{2i}=\lambda_i(1+\delta_i).
\]

\emph{Factorial moments.}
The product $\prod_i(X_i)_{r_i}$ counts ordered tuples of $r=\sum_ir_i$ pairwise distinct cycles, with $r_i$ cycles of length $i$. It equals $\prod_i(2i)^{-r_i}$ times the number of tuples of directed closed walks with pairwise distinct underlying cycles.

First consider tuples whose cycles are pairwise vertex-disjoint. Their vertices can be chosen in $(n)_{u}=n^{u}(1+O(n^{-1}))$ ways, where $u=\sum_i i\,r_i$. The $u$ prescribed pairs are pairwise disjoint, so~\eqref{eq:pair-prob} applies with $a$ and $b$ equal to the total numbers of $F$- and $R$-pairs. The weight therefore factorises over the cycles, and these tuples contribute $\prod_i(\tr T^i)^{r_i}(1+O(n^{-1}))$. Dividing by $\prod_i(2i)^{r_i}$ gives the right-hand side of~\eqref{eq:factorial-moments}.

Now consider a tuple in which two cycles share a vertex. Let $U$ be the union of its cycles, a multigraph with $u'\le u$ vertices and $e'$ edges, where the edges are distinct pairs of $P$. There are $O(1)$ possible isomorphism types, $O(n^{u'})$ choices of vertex labels and points for each type, and probability $O(n^{-e'})$ that the $e'$ prescribed pairs belong to $P$, by~\eqref{eq:pair-prob}. This probability is zero unless the pairs are disjoint. Every connected component $C$ of $U$ contains a cycle, so $e(C)\ge u(C)$. If $e(C)=u(C)$, then $C$ has cyclomatic number $e(C)-u(C)+1=1$ and contains exactly one cycle, with loops and pairs of parallel edges counted as cycles. Some component of $U$ contains two distinct cycles from the tuple, so that component has $e(C)\ge u(C)+1$. Therefore $e'\ge u'+1$, and tuples with a shared vertex contribute $O(n^{u'-e'})=O(n^{-1})$. This proves~\eqref{eq:factorial-moments}.
\end{proof}

\section{The second moment}\label{sec:second-moment}

In this section, we fix $d\ge3$ and $1\le k\le d-1$ and prove the following lemma.

\begin{lemma}[second moment]\label{lem:second-moment}
As $n\to\infty$ with $dn$ and $kn$ even,
\begin{equation}\label{eq:second-moment-pairing}
\frac{\E\bigl[Z_k(P)^2\bigr]}{\bigl(\E Z_k(P)\bigr)^2}
\;\longrightarrow\;
\sqrt{\frac{d-1}{d-2}}
\;=\;
\exp\Bigl(\sum_{i=1}^{\infty}\lambda_i\delta_i^2\Bigr).
\end{equation}
\end{lemma}

For every pairing $P$, complementation $F\mapsto P\setminus F$ is a bijection between the $k$-factors and the $(d-k)$-factors of $P$. Thus
\begin{equation}\label{eq:complement}
Z_{d-k}(P)=Z_k(P)\qquad\text{for every }P\in\mathcal P_{n,d}\text{ and every }n .
\end{equation}
Both sides of~\eqref{eq:second-moment-pairing} are invariant under $k\mapsto d-k$, so it suffices to prove Lemma~\ref{lem:second-moment} for $k\le d/2$. \emph{Throughout this section, we assume that $1\le k\le d/2$}; since $d\ge3$, this excludes $(d,k)=(2,1)$. This reduction requires no asymptotic input such as contiguity: \eqref{eq:complement} is an identity of random variables, so all moments and limit distributions of $Z_k(P)$ and $Z_{d-k}(P)$ coincide exactly.

Let $\Hyp(d,k,k)$ denote the hypergeometric distribution of the overlap of two independent uniform $k$-subsets of a $d$-set:
\[
\Prob(H=j)=\frac{\binom{k}{j}\binom{d-k}{k-j}}{\binom{d}{k}}\qquad(0\le j\le k),
\]
and all $k+1$ values have positive probability because $k\le d-k$. Its mean and variance are
\begin{equation}\label{eq:hyp-moments}
\mu_H:=\E H=\frac{k^2}{d},\qquad
\sigma_H^2:=\Var H=\frac{k^2(d-k)^2}{d^2(d-1)} .
\end{equation}

\subsection{Reduction to a lattice sum}

\begin{lemma}\label{lem:reduction}
Let $H_1,\dots,H_n$ be i.i.d.\ $\Hyp(d,k,k)$ and $S_n=\sum_{v=1}^nH_v$. Then
\begin{equation}\label{eq:exact-second}
\frac{\E\bigl[Z_k(P)^2\bigr]}{\bigl(\E Z_k(P)\bigr)^2}
=\sum_{\substack{0\le\sigma\le kn \\ \sigma\text{ even}}}
\Prob(S_n=\sigma)\,R_M(\sigma),
\quad
R_M(\sigma):=\frac{M(\sigma)\,M(kn-\sigma)^2\,M\bigl((d-2k)n+\sigma\bigr)\,M(dn)}{M(kn)^2\,M\bigl((d-k)n\bigr)^2}.
\end{equation}
\end{lemma}

\begin{proof}
$Z_k(P)^2$ counts ordered pairs $(F_1,F_2)$ of $k$-factors of $P$. Let $A_v,B_v\subseteq V_v$ be the $k$-subsets of points covered by $F_1$ and $F_2$, and put $\mathcal A=\bigcup_vA_v$ and $\mathcal B=\bigcup_vB_v$. Given the sets $(A_v,B_v)_v$, the factors are determined by $P$: $F_1$ is the set of pairs of $P$ within $\mathcal A$, and $F_2$ is the set of pairs within $\mathcal B$. Conversely, $P$ \emph{extends} $(A_v,B_v)_v$ to a pair of $k$-factors if and only if $P$ pairs $\mathcal A$ internally and pairs $\mathcal B$ internally.

We claim that this condition holds if and only if $P$ pairs each of the four classes
\[
\mathcal A\cap\mathcal B,\qquad \mathcal A\setminus\mathcal B,\qquad \mathcal B\setminus\mathcal A,\qquad (\mathcal A\cup\mathcal B)^{c}
\]
internally. One direction is immediate. Conversely, suppose that $\mathcal A$ and $\mathcal B$ are paired internally. If $x\in\mathcal A\cap\mathcal B$, then its partner belongs to both $\mathcal A$ and $\mathcal B$, and hence to $\mathcal A\cap\mathcal B$. If $x\in\mathcal A\setminus\mathcal B$, then its partner $x'$ lies in $\mathcal A$. If $x'\in\mathcal B$, then $x'\in\mathcal A\cap\mathcal B$, so its partner $x$ would lie in $\mathcal A\cap\mathcal B$, a contradiction. Hence $x'\in\mathcal A\setminus\mathcal B$. The argument for $\mathcal B\setminus\mathcal A$ is symmetric, and the remaining class is then also paired internally.

Set $H_v:=|A_v\cap B_v|$ and $\sigma:=\sum_vH_v$. The four classes have sizes $\sigma$, $kn-\sigma$, $kn-\sigma$, and $(d-2k)n+\sigma$. The number of pairings that extend $(A_v,B_v)_v$ is therefore
\[
M(\sigma)\,M(kn-\sigma)^2\,M\bigl((d-2k)n+\sigma\bigr).
\]
Since $kn$ and $(d-2k)n=dn-2kn$ are even, all four sizes are even exactly when $\sigma$ is even. Otherwise, there are no extensions. Summing over the $\binom dk^{2n}$ choices of $(A_v,B_v)_v$ and dividing by $M(dn)$ gives
\[
\E\bigl[Z_k(P)^2\bigr]=\frac{\binom{d}{k}^{2n}}{M(dn)}\,\E\Bigl[M(S_n)\,M(kn-S_n)^2\,M\bigl((d-2k)n+S_n\bigr)\Bigr],
\]
where the overlaps $H_v=|A_v\cap B_v|$ of independent uniform $k$-subsets $A_v,B_v\subseteq V_v$ are i.i.d.\ $\Hyp(d,k,k)$. Dividing by $(\E Z_k(P))^2=\binom dk^{2n}M(kn)^2M((d-k)n)^2/M(dn)^2$, from~\eqref{eq:EZk-pairing}, gives~\eqref{eq:exact-second}.
\end{proof}

\subsection{Stirling bounds}

Define, with the convention $0\log0=0$,
\begin{equation}\label{eq:g-def}
g(s)\;:=\;\frac{s}{2}\log s+(k-s)\log(k-s)+\frac{d-2k+s}{2}\log(d-2k+s)+\frac{d}{2}\log d-k\log k-(d-k)\log(d-k)
\end{equation}
for $0\le s\le k$.

\begin{lemma}\label{lem:stirling}
\begin{enumerate}
\item[\textup{(a)}] For every even $\sigma$ with $0\le\sigma\le kn$, $R_M(\sigma)\le 32\,e^{ng(\sigma/n)}$.
\item[\textup{(b)}] For every compact $J\subset(0,k)$, $R_M(\sigma)=\sqrt2\,e^{ng(\sigma/n)}\bigl(1+O(n^{-1})\bigr)$ uniformly over even $\sigma$ with $\sigma/n\in J$.
\item[\textup{(c)}] $g$ is continuous on $[0,k]$ and smooth on $(0,k)$, with $g(\mu_H)=g'(\mu_H)=0$, 
\[
g''(s)=\frac{1}{2s}+\frac{1}{k-s}+\frac{1}{2(d-2k+s)}\quad(0<s<k),\qquad
g''(\mu_H)=\frac{d^3}{2k^2(d-k)^2}.
\]
\end{enumerate}
\end{lemma}

\begin{proof}
(a) The bounds $\sqrt{2\pi N}(N/e)^N\le N!\le e\sqrt N(N/e)^N$ for $N\ge1$ give, for all $\ell\ge1$,
\[
\frac{2\sqrt\pi}{e}\Bigl(\frac{2\ell}{e}\Bigr)^{\ell}\le M(2\ell)=\frac{(2\ell)!}{\ell!\,2^{\ell}}\le\frac{e}{\sqrt\pi}\Bigl(\frac{2\ell}{e}\Bigr)^{\ell},
\qquad\text{hence}\qquad
\Bigl(\frac{2\ell}{e}\Bigr)^{\ell}\le M(2\ell)\le 2\Bigl(\frac{2\ell}{e}\Bigr)^{\ell}
\]
for all $\ell\ge0$, with $0^0=1$. For $\ell=0$, the last display reads $1\le1\le2$. Apply the upper bound to the five factors in the numerator of $R_M(\sigma)$ and the lower bound to the four factors in the denominator. The exponents of $e$ and of $n$ cancel because
\[
\frac{\sigma}{2}+(kn-\sigma)+\frac{(d-2k)n+\sigma}{2}+\frac{dn}{2}=dn=kn+(d-k)n,
\]
and hence $R_M(\sigma)\le 2^5e^{ng(\sigma/n)}$. Here $g$ collects the logarithms of $\sigma/n$, $(kn-\sigma)/n$, $((d-2k)n+\sigma)/n$, $d$, $k$, and $d-k$, with the coefficients in~\eqref{eq:g-def}. Some of the arguments $\sigma$, $kn-\sigma$, and $(d-2k)n+\sigma$ may vanish; the case $\ell=0$ and the convention $0\log0=0$ cover these possibilities.

(b) If $\sigma/n\in J$, all nine arguments are of order $n$. Applying $M(2\ell)=\sqrt2\,(2\ell/e)^{\ell}(1+O(\ell^{-1}))$ to each factor gives the result, with the factor $(\sqrt2)^{5-4}=\sqrt2$.

(c) Continuity on $[0,k]$ follows from $x\log x\to0$ as $x\downarrow0$. Differentiation gives
\[
g'(s)=\tfrac12\log s-\log(k-s)+\tfrac12\log(d-2k+s),
\]
and $g''$ has the stated form. At $\mu_H=k^2/d$, we have $k-\mu_H=k(d-k)/d$ and $d-2k+\mu_H=(d-k)^2/d$. Thus $g'(\mu_H)=\log k-\tfrac12\log d-\log k-\log(d-k)+\log d+\log(d-k)-\tfrac12\log d=0$. In $g(\mu_H)$, the coefficients of $\log k$, $\log(d-k)$, and $\log d$ are $\frac{k^2}{d}+\frac{k(d-k)}{d}-k=0$, $\frac{k(d-k)}{d}+\frac{(d-k)^2}{d}-(d-k)=0$, and $-\frac{k^2+2k(d-k)+(d-k)^2}{2d}+\frac d2=0$, respectively. Finally,
\[
g''(\mu_H)=\frac{d}{2k^2}+\frac{d}{k(d-k)}+\frac{d}{2(d-k)^2}=\frac d2\Bigl(\frac1k+\frac1{d-k}\Bigr)^2=\frac{d^3}{2k^2(d-k)^2}.
\qedhere
\]
\end{proof}

\subsection{The hypergeometric distribution is a Poisson-binomial distribution}

\begin{lemma}\label{lem:poisson-binomial}
Let $1\le k\le d/2$ and $H\sim\Hyp(d,k,k)$. There exist distinct $p_1,\dots,p_k\in(0,1)$ such that
\[
\E\bigl[t^{H}\bigr]=\prod_{i=1}^k\bigl(1-p_i+p_it\bigr)\qquad(t\in\mathbb R),
\]
that is, $H$ has the distribution of a sum of $k$ independent $\mathrm{Bernoulli}(p_i)$ variables. In particular, $\sum_ip_i=\mu_H$ and $\sum_ip_i(1-p_i)=\sigma_H^2$.
\end{lemma}

\begin{proof}
Consider the polynomials
\[
p(u):=u^{d-k}(u-1)^k,
\qquad
Q(u):=\sum_{j=0}^k\binom kj\binom{d-k}{k-j}u^j(u-1)^{k-j}.
\]
The Leibniz rule gives
\[
\frac{1}{k!}\,\frac{d^k}{du^k}\bigl[u^{d-k}(u-1)^k\bigr]
=\sum_{i=0}^k\binom{k}{i}\binom{d-k}{i}u^{d-k-i}(u-1)^{i}
=u^{d-2k}\sum_{j=0}^k\binom{k}{j}\binom{d-k}{k-j}u^{j}(u-1)^{k-j}
\]
where we substituted $j=k-i$. Hence $\frac{1}{k!}p^{(k)}(u)=u^{d-2k}Q(u)$.

We claim that, for $0\le i\le k$, the polynomial $p^{(i)}$, of degree $d-i$, has a zero of multiplicity exactly $d-k-i$ at $0$, a zero of multiplicity exactly $k-i$ at $1$, and $i$ distinct simple zeros in $(0,1)$. These are all its zeros. The claim holds for $i=0$. Suppose it holds for some $i<k$. Then $p^{(i+1)}$ vanishes to order exactly $d-k-i-1\ge0$ at $0$ and $k-i-1\ge0$ at $1$. By Rolle's theorem, it has a zero strictly between every two consecutive zeros of $p^{(i)}$ in $[0,1]$. This gives $i+1$ zeros in $(0,1)$: between $0$ and the smallest interior zero, between consecutive interior zeros, and between the largest interior zero and $1$. The listed zeros number $(d-k-i-1)+(k-i-1)+(i+1)=d-i-1=\deg p^{(i+1)}$. They are therefore all the zeros, and the interior zeros are simple.

For $i=k$, we conclude that $u^{d-2k}Q(u)$ has exactly $k$ simple zeros $u_1<\dots<u_k$ in $(0,1)$ besides $u=0$. Thus $Q$ has the $k$ simple zeros $u_1,\dots,u_k\in(0,1)$. Since $\deg Q\le k$ and the coefficient of $u^k$ in $Q$ is $\sum_j\binom kj\binom{d-k}{k-j}=\binom dk$,
\[
Q(u)=\binom dk\prod_{i=1}^k(u-u_i).
\]
On the other hand, for $u\ne1$,
\[
Q(u)=(u-1)^k\sum_{j=0}^k\binom kj\binom{d-k}{k-j}\Bigl(\frac{u}{u-1}\Bigr)^{j}
=\binom dk(u-1)^k\,\E\Bigl[\Bigl(\frac{u}{u-1}\Bigr)^{H}\Bigr].
\]
Compare the two expressions and substitute $t=u/(u-1)$, so $u-1=1/(t-1)$ and $(u-u_i)/(u-1)=(1-u_i)t+u_i$. We obtain
\[
\E\bigl[t^{H}\bigr]=\prod_{i=1}^k\frac{u-u_i}{u-1}=\prod_{i=1}^k\bigl(u_i+(1-u_i)t\bigr)
\]
for all $t\ne1$. At $t=1$ both sides equal $1$, so the identity holds for all $t$. Set $p_i:=1-u_i\in(0,1)$. The $p_i$ are distinct, and
\[
\E\bigl[t^{H}\bigr]=\prod_{i=1}^k\bigl(1-p_i+p_it\bigr)
\]
for all $t$. Thus $H$ has the distribution of a sum of $k$ independent Bernoulli($p_i$) variables, with mean $\sum_i p_i$ and variance $\sum_i p_i(1-p_i)$. Hence $\sum_i p_i=\mu_H$ and $\sum_i p_i(1-p_i)=\sigma_H^2$, by~\eqref{eq:hyp-moments}.
\end{proof}

\subsection{Exponential tilting and the rate function}

Let $K(\theta):=\log\E e^{\theta H}$, $\theta\in\mathbb R$, be the cumulant generating function of $H\sim\Hyp(d,k,k)$. For $\theta\in\mathbb R$, let $\Prob_\theta$ be the probability measure under which $H_1,\dots,H_n$ are i.i.d.\ with the \emph{tilted} distribution
\[
\Prob_\theta(H_1=j)=\Prob(H_1=j)\,e^{\theta j-K(\theta)}\qquad(0\le j\le k).
\]
By Lemma~\ref{lem:poisson-binomial}, $K(\theta)=\sum_{i=1}^k\log(1-p_i+p_ie^\theta)$, and the tilted probability generating function is
\begin{equation}\label{eq:tilted-pgf}
\E_\theta\bigl[t^{H_1}\bigr]=\frac{\E[(te^\theta)^{H_1}]}{\E[e^{\theta H_1}]}=\prod_{i=1}^k\bigl(1-q_i(\theta)+q_i(\theta)t\bigr),
\qquad
q_i(\theta):=\frac{p_ie^\theta}{1-p_i+p_ie^\theta},
\end{equation}
so that, under $\Prob_\theta$, $H_1$ is again a sum of $k$ independent Bernoulli variables. Their parameters $q_i(\theta)$ are distinct for each fixed $\theta$ because $q_i(\theta)$ is strictly increasing in $p_i$. Consequently,
\begin{equation}\label{eq:K-derivatives}
K'(\theta)=\E_\theta H_1=\sum_{i=1}^kq_i(\theta),\qquad
K''(\theta)=\Var_\theta H_1=\sum_{i=1}^kq_i(\theta)\bigl(1-q_i(\theta)\bigr)>0 .
\end{equation}
Since $q_i(\theta)\to0$ as $\theta\to-\infty$ and $q_i(\theta)\to1$ as $\theta\to\infty$, $K'$ is a smooth increasing bijection from $\mathbb R$ to $(0,k)$. Let $\theta(s)$, $0<s<k$, be its inverse. This function is smooth, and $\theta(\mu_H)=0$ because $K'(0)=\E H=\mu_H$. Define the rate function and the comparison function by
\begin{equation}\label{eq:I-Phi-def}
\rho(s):=\sup_{\theta\in\mathbb R}\bigl(\theta s-K(\theta)\bigr)\quad(0\le s\le k),
\qquad
\Phi(s):=g(s)-\rho(s).
\end{equation}
For $0<s<k$, the supremum is attained at $\theta=\theta(s)$ because the derivative of the concave function $\theta\mapsto\theta s-K(\theta)$ is $s-K'(\theta)$. Thus $\rho(s)=\theta(s)s-K(\theta(s))$.

\begin{lemma}\label{lem:rate}
\begin{enumerate}
\item[\textup{(a)}] $\rho$ is smooth on $(0,k)$ with $\rho'(s)=\theta(s)$ and $\rho''(s)=1/K''(\theta(s))$, and $\rho(\mu_H)=\rho'(\mu_H)=0$.
\item[\textup{(b)}] For $0<s<k$,
\[
\rho''(s)=\frac{1}{K''(\theta(s))}\;\ge\;\frac1s+\frac1{k-s}\;\ge\;g''(s),
\qquad
\frac1s+\frac1{k-s}-g''(s)=\frac{d-2k}{2s(d-2k+s)},
\]
and at least one of the two inequalities is strict. Hence $\Phi''<0$ on $(0,k)$, $\Phi(\mu_H)=\Phi'(\mu_H)=0$, and
\[
\Phi''(\mu_H)=g''(\mu_H)-\frac{1}{\sigma_H^2}=-\frac{d^2(d-2)}{2k^2(d-k)^2}<0 .
\]
\item[\textup{(c)}] $\Phi$ is continuous on $[0,k]$, and there is $c=c(d,k)>0$ such that $\Phi(s)\le-c\,(s-\mu_H)^2$ for all $0\le s\le k$.
\item[\textup{(d)}] For every $\theta\in\mathbb R$ and every integer $\sigma$ with $0\le\sigma\le kn$, writing $s=\sigma/n$,
\begin{equation}\label{eq:tilting-identity}
\Prob(S_n=\sigma)=e^{-n(\theta s-K(\theta))}\,\Prob_\theta(S_n=\sigma),
\qquad\text{hence}\qquad
\Prob(S_n=\sigma)\le e^{-n\rho(s)}
\end{equation}
and, for even $\sigma$, $\Prob(S_n=\sigma)R_M(\sigma)\le 32\,e^{n\Phi(s)}$.
\end{enumerate}
\end{lemma}

\begin{proof}
(a) Differentiating $\rho(s)=\theta(s)s-K(\theta(s))$ gives $\rho'(s)=\theta(s)+\theta'(s)\bigl(s-K'(\theta(s))\bigr)=\theta(s)$. The inverse function theorem gives $\rho''(s)=\theta'(s)=1/K''(\theta(s))$. At $\mu_H$, $\theta(\mu_H)=0$ and $K(0)=0$, so $\rho(\mu_H)=\rho'(\mu_H)=0$.

(b) Put $q_i=q_i(\theta(s))$, so $\sum_iq_i=s$ by~\eqref{eq:K-derivatives}. The Cauchy--Schwarz inequality gives $\sum_iq_i^2\ge(\sum_iq_i)^2/k=s^2/k$, with equality if and only if all $q_i$ are equal. Hence
\[
K''(\theta(s))=\sum_iq_i(1-q_i)=s-\sum_iq_i^2\le s-\frac{s^2}{k}=\frac{s(k-s)}{k},
\qquad\text{i.e.}\qquad
\frac{1}{K''(\theta(s))}\ge\frac1s+\frac1{k-s}.
\]
If $k\ge2$, the $q_i$ are distinct, so this inequality is strict. The identity for $\frac1s+\frac1{k-s}-g''(s)$ follows from Lemma~\ref{lem:stirling}(c): $\frac1s-\frac1{2s}-\frac1{2(d-2k+s)}=\frac{(d-2k+s)-s}{2s(d-2k+s)}$. It is nonnegative and is strictly positive if $d>2k$. If $k=1$, then $d\ge3>2k$. Thus $\rho''(s)>g''(s)$ for $0<s<k$, that is, $\Phi''(s)<0$.

Part~(a) and Lemma~\ref{lem:stirling}(c) give $\Phi(\mu_H)=g(\mu_H)-\rho(\mu_H)=0$ and $\Phi'(\mu_H)=g'(\mu_H)-\theta(\mu_H)=0$. Finally, $K''(0)=\Var H=\sigma_H^2$, so
\[
\Phi''(\mu_H)=g''(\mu_H)-\frac1{\sigma_H^2}=\frac{d^3}{2k^2(d-k)^2}-\frac{d^2(d-1)}{k^2(d-k)^2}=-\frac{d^2(d-2)}{2k^2(d-k)^2}.
\]

(c) The function $g$ is continuous on $[0,k]$ by Lemma~\ref{lem:stirling}(c). For $\rho$, since $K$ is increasing, $\rho(0)=-\inf_\theta K(\theta)=-\lim_{\theta\to-\infty}K(\theta)=-\log\Prob(H=0)=\log\bigl(\binom dk/\binom{d-k}k\bigr)$, which is finite because $k\le d-k$. Similarly, $\rho(k)=-\log\Prob(H=k)=\log\binom dk$. The function $\rho$ is a finite convex function on $[0,k]$ because it is a supremum of affine functions. It is therefore continuous on $(0,k)$ and lower semicontinuous on $[0,k]$. Convexity also gives $\rho(s)\le\frac{k-s}{k}\rho(0)+\frac sk\rho(k)$, and hence $\limsup_{s\downarrow0}\rho(s)\le \rho(0)$ and $\limsup_{s\uparrow k}\rho(s)\le \rho(k)$. Thus $\rho$, and hence $\Phi$, is continuous on $[0,k]$.

By part~(b), $\Phi$ is strictly concave on $(0,k)$, with $\Phi(\mu_H)=\Phi'(\mu_H)=0$. Thus $\Phi$ is strictly increasing on $(0,\mu_H]$ and strictly decreasing on $[\mu_H,k)$. By continuity, $\Phi(s)<0$ for all $s\in[0,k]\setminus\{\mu_H\}$. Taylor's theorem and $\Phi''(\mu_H)<0$ give $\kappa>0$ such that $\Phi(s)\le\frac14\Phi''(\mu_H)(s-\mu_H)^2$ for $|s-\mu_H|\le\kappa$. By compactness, $\eta:=-\max\{\Phi(s):s\in[0,k],\,|s-\mu_H|\ge\kappa\}>0$. The claim follows with $c=\min\{|\Phi''(\mu_H)|/4,\;\eta/k^2\}$.

(d) The density of $\Prob$ with respect to $\Prob_\theta$ on $(H_1,\dots,H_n)$ is $e^{-\theta S_n+nK(\theta)}$. Hence $\Prob(S_n=\sigma)=\E_\theta\bigl[e^{-\theta S_n+nK(\theta)}\mathbf 1\{S_n=\sigma\}\bigr]=e^{-\theta\sigma+nK(\theta)}\Prob_\theta(S_n=\sigma)$. Bounding $\Prob_\theta(S_n=\sigma)\le1$ and taking the supremum over $\theta$ gives $\Prob(S_n=\sigma)\le e^{-n\rho(s)}$. Combining this bound with Lemma~\ref{lem:stirling}(a) proves the last assertion.
\end{proof}

\subsection{A uniform local limit theorem for the tilted sums}

\begin{lemma}\label{lem:llt}
Fix $\theta_0>0$. There is a sequence $\varepsilon_n\to0$, depending only on $d,k,\theta_0$, such that for all $\theta\in[-\theta_0,\theta_0]$ and all $n$ for which $\sigma:=nK'(\theta)$ is an integer,
\[
\Bigl|\sqrt{2\pi nK''(\theta)}\;\Prob_\theta(S_n=\sigma)-1\Bigr|\le\varepsilon_n .
\]
\end{lemma}

\begin{proof}
Fix $\theta\in[-\theta_0,\theta_0]$ and write $\mu=K'(\theta)$ and $v=K''(\theta)$. Then $\sigma=n\mu$ and $v\in[v_0,k^2]$, where $v_0:=\min_{|\theta|\le\theta_0}K''(\theta)>0$ by~\eqref{eq:K-derivatives} and continuity. Let $\phi_\theta(t)=\E_\theta e^{itH_1}$ and $\psi_\theta(t)=e^{-it\mu}\phi_\theta(t)=\E_\theta e^{it(H_1-\mu)}$. By~\eqref{eq:tilted-pgf}, $\phi_\theta(t)=\prod_i(1-q_i+q_ie^{it})$, where $q_i=q_i(\theta)$. Since $|1-q+qe^{it}|^2=1-2q(1-q)(1-\cos t)\le e^{-2q(1-q)(1-\cos t)}$,
\begin{equation}\label{eq:cf-bound}
|\psi_\theta(t)|=|\phi_\theta(t)|\le e^{-v(1-\cos t)}\le e^{-2vt^2/\pi^2}\qquad(|t|\le\pi),
\end{equation}
where the last inequality uses $1-\cos t\ge2t^2/\pi^2$ on $[-\pi,\pi]$.

Since $|H_1-\mu|\le k$, the inequality $|e^{ix}-1-ix+x^2/2|\le|x|^3/6$ gives $\psi_\theta(t)=1+z(t)$ with $|z(t)+vt^2/2|\le k^3|t|^3/6$. Choose $t_0\in(0,\pi]$, depending only on $k$, such that $k^2t_0^2/2+k^3t_0^3/6\le1/2$. Then $|z(t)|\le\tfrac12$ for $|t|\le t_0$. Using $|\log(1+z)-z|\le|z|^2$ for $|z|\le\tfrac12$, we obtain a constant $C=C(k)$ such that
\begin{equation}\label{eq:log-psi}
\Bigl|\log\psi_\theta(t)+\frac{vt^2}{2}\Bigr|\le C|t|^3\qquad(|t|\le t_0),
\end{equation}
where $\log$ is the principal branch.

Fourier inversion for the integer-valued variable $S_n$, followed by the substitution $t=\tau/\sqrt{nv}$, gives
\[
\sqrt{2\pi nv}\;\Prob_\theta(S_n=n\mu)
=\sqrt{2\pi nv}\,\frac{1}{2\pi}\int_{-\pi}^{\pi}e^{-it n\mu}\phi_\theta(t)^n\,dt
=\frac{1}{\sqrt{2\pi}}\int_{|\tau|\le\pi\sqrt{nv}}\psi_\theta\bigl(\tau/\sqrt{nv}\bigr)^n\,d\tau .
\]
Split the integral at $|\tau|=n^{1/8}$. For $|\tau|\le n^{1/8}$ and $n$ large enough that $n^{1/8}/\sqrt{nv_0}\le t_0$, \eqref{eq:log-psi} gives
\[
\Bigl|n\log\psi_\theta\bigl(\tau/\sqrt{nv}\bigr)+\frac{\tau^2}{2}\Bigr|\le\frac{Cn|\tau|^3}{(nv)^{3/2}}\le Cv_0^{-3/2}n^{-1/8}=:\varepsilon_n',
\]
so
\[
\Bigl|\psi_\theta\bigl(\tau/\sqrt{nv}\bigr)^n-e^{-\tau^2/2}\Bigr|\le e^{-\tau^2/2}\bigl(e^{\varepsilon_n'}-1\bigr).
\]
Hence the integral over $|\tau|\le n^{1/8}$ equals $\int_{|\tau|\le n^{1/8}}e^{-\tau^2/2}d\tau+O(\varepsilon_n')=\sqrt{2\pi}+O(\varepsilon_n')+O(e^{-n^{1/4}/2})$. For $n^{1/8}<|\tau|\le\pi\sqrt{nv}$, \eqref{eq:cf-bound} gives $|\psi_\theta(\tau/\sqrt{nv})|^n\le e^{-2\tau^2/\pi^2}$, so this part of the integral is $O(e^{-2n^{1/4}/\pi^2})$. All constants depend only on $d,k,\theta_0$. This proves the lemma.
\end{proof}

\subsection{Proof of Lemma~\ref{lem:second-moment}}

\begin{proof}[Proof of Lemma~\ref{lem:second-moment}]
Recall that $1\le k\le d/2$. Let $a:=\frac{1}{\sigma_H^2}-g''(\mu_H)=-\Phi''(\mu_H)=\frac{d^2(d-2)}{2k^2(d-k)^2}>0$ by Lemma~\ref{lem:rate}(b). Split the sum in~\eqref{eq:exact-second} into $\Sigma_{\mathrm{out}}$, over even $\sigma$ with $|\sigma-n\mu_H|>n^{3/5}$, and $\Sigma_{\mathrm{in}}$, over even $\sigma$ with $|\sigma-n\mu_H|\le n^{3/5}$.

\emph{The tail.}
By Lemma~\ref{lem:rate}(c),(d), every term of $\Sigma_{\mathrm{out}}$ is at most $32\,e^{n\Phi(\sigma/n)}\le32\,e^{-cn(\sigma/n-\mu_H)^2}\le32\,e^{-cn^{1/5}}$. There are at most $kn+1$ terms, so $\Sigma_{\mathrm{out}}=O(ne^{-cn^{1/5}})\to0$.

\emph{The central window.}
Let $\sigma$ be even with $|\sigma-n\mu_H|\le n^{3/5}$, and put $s=\sigma/n$ and $x=(\sigma-n\mu_H)/\sqrt n$. Then $|s-\mu_H|\le n^{-2/5}$ and $|x|\le n^{1/10}$. For large $n$, $s$ lies in a fixed compact subset of $(0,k)$, and $\theta(s)=O(n^{-2/5})$ lies in $[-1,1]$. Apply the tilting identity~\eqref{eq:tilting-identity} with $\theta=\theta(s)$, for which $nK'(\theta(s))=\sigma$ and $\theta(s)s-K(\theta(s))=\rho(s)$. Lemma~\ref{lem:llt}, with $\theta_0=1$, then gives
\[
\Prob(S_n=\sigma)=e^{-n\rho(s)}\,\frac{1+o(1)}{\sqrt{2\pi nK''(\theta(s))}}=e^{-n\rho(s)}\,\frac{1+o(1)}{\sqrt{2\pi n\sigma_H^2}}
\]
uniformly in the window, because $K''(\theta(s))=K''(0)+O(n^{-2/5})=\sigma_H^2+O(n^{-2/5})$. Lemma~\ref{lem:stirling}(b) gives $R_M(\sigma)=\sqrt2\,e^{ng(s)}(1+O(n^{-1}))$. Therefore, uniformly in the window,
\[
\Prob(S_n=\sigma)\,R_M(\sigma)=\frac{\sqrt2}{\sqrt{2\pi n\sigma_H^2}}\,e^{n\Phi(s)}\bigl(1+o(1)\bigr)=\frac{2}{\sqrt{4\pi n\sigma_H^2}}\,e^{n\Phi(s)}\bigl(1+o(1)\bigr).
\]
The function $\Phi$ is smooth near $\mu_H$, with $\Phi(\mu_H)=\Phi'(\mu_H)=0$ and $\Phi''(\mu_H)=-a$. Taylor's theorem gives
\[
n\Phi(s)=-\frac{a}{2}\,n(s-\mu_H)^2+O\bigl(n|s-\mu_H|^3\bigr)=-\frac{a}{2}\,x^2+O\bigl(n^{-1/5}\bigr).
\]
Consequently,
\[
\Sigma_{\mathrm{in}}=\bigl(1+o(1)\bigr)\,\frac{1}{\sqrt{4\pi\sigma_H^2}}\sum_{\substack{|\sigma-n\mu_H|\le n^{3/5}\\ \sigma\text{ even}}}\frac{2}{\sqrt n}\,e^{-ax^2/2}.
\]
Consecutive even values of $\sigma$ correspond to values of $x$ spaced by $2/\sqrt n$, and the window corresponds to $|x|\le n^{1/10}$. For fixed $A>0$, the terms with $|x|\le A$ form a Riemann sum that converges to $\int_{-A}^{A}e^{-ax^2/2}\,dx$. Since $e^{-ax^2/2}$ is monotone on each half-line, the terms with $|x|>A$ sum to at most $\int_{|x|>A-2/\sqrt n}e^{-ax^2/2}\,dx$. Letting $n\to\infty$ and then $A\to\infty$, we see that $\Sigma_{\mathrm{in}}$ converges to
\[
\frac{1}{\sqrt{4\pi\sigma_H^2}}\int_{-\infty}^{\infty}e^{-ax^2/2}\,dx=\frac{1}{\sqrt{4\pi\sigma_H^2}}\sqrt{\frac{2\pi}{a}}=\frac{1}{\sqrt{2a\sigma_H^2}} .
\]
Finally, by~\eqref{eq:hyp-moments} and Lemma~\ref{lem:stirling}(c), $2\sigma_H^2g''(\mu_H)=\frac{2k^2(d-k)^2}{d^2(d-1)}\cdot\frac{d^3}{2k^2(d-k)^2}=\frac{d}{d-1}$, so
\[
2a\sigma_H^2=2-2\sigma_H^2g''(\mu_H)=2-\frac{d}{d-1}=\frac{d-2}{d-1},
\qquad
\lim_{n\to\infty}\frac{\E[Z_k(P)^2]}{(\E Z_k(P))^2}=\sqrt{\frac{d-1}{d-2}} .
\]
The identity with $\exp(\sum_{i\ge1}\lambda_i\delta_i^2)$ was noted after~\eqref{eq:lambda-delta}.
\end{proof}

\begin{remark}
Within the window, $\Prob(S_n=\sigma)$ can be as small as $e^{-cn^{1/5}}$, far below the $n^{-1/2}$ scale of the ordinary local central limit theorem. Its additive error is therefore too large. The tilting identity~\eqref{eq:tilting-identity} reduces the problem to a local limit theorem \emph{at the mean} of the tilted distribution. Lemma~\ref{lem:llt} then gives a relative error $o(1)$ uniformly in the tilt. The strict inequality $\Phi''<0$ in Lemma~\ref{lem:rate}(b) relies on the Poisson-binomial representation. When $d=2k$, so that $g''(s)=\frac1s+\frac1{k-s}$ exactly, strictness follows from the distinct Bernoulli parameters.
\end{remark}

\section{Proof of Theorem~\ref{thm:main}}\label{sec:proof}

Throughout this section, $d\ge3$ and $1\le k\le d-1$ are fixed, and $n\to\infty$ along integers for which $dn$ and $kn$ are even. Let
\[
Y_{k,n}:=\frac{Z_k(P)}{\E Z_k(P)},\qquad P\sim\mathcal P_{n,d},
\]
where $X_i=X_i(P)$ is defined as in Section~\ref{sec:pairing}. Write $W^{\mathrm{pair}}:=\prod_{i\ge1}(1+\delta_i)^{\xi_i}e^{-\lambda_i\delta_i}$, $W:=\prod_{i\ge3}(1+\delta_i)^{\xi_i}e^{-\lambda_i\delta_i}$ as in~\eqref{eq:W-def}, and $W^{(m)}(\xi):=\prod_{i\le m}(1+\delta_i)^{\xi_i}e^{-\lambda_i\delta_i}$. By Lemma~\ref{lem:W-properties}, these products converge almost surely and in $L^2$, are almost surely positive, and have mean $1$. Moreover, $W$ is independent of $(\xi_1,\xi_2)$.

\medskip\noindent
\emph{Step 1: small subgraph conditioning in the pairing model.}
Hypothesis (i) of Theorem~\ref{thm:janson} is~\eqref{eq:poisson-cycles}, and (ii) is Lemma~\ref{lem:short-cycles}. Hypothesis (iii) holds because $\sum_i\lambda_i\delta_i^2=\frac12\log\frac{d-1}{d-2}$, and (iv) is Lemma~\ref{lem:second-moment}, which applies to all $1\le k\le d-1$ by~\eqref{eq:complement}. Also, $\delta_i>-1$ for all $i$. Hence, for every fixed $m$,
\begin{equation}\label{eq:joint-pairing}
\bigl(Y_{k,n},X_1,\dots,X_m\bigr)\;\Longrightarrow\;\bigl(W^{\mathrm{pair}},\xi_1,\dots,\xi_m\bigr)\qquad\text{in }\mathcal P_{n,d}.
\end{equation}

\medskip\noindent
\emph{Step 2: uniform integrability.}
By~\eqref{eq:joint-pairing}, $Y_{k,n}^2\Rightarrow(W^{\mathrm{pair}})^2$. By (iv) and Lemma~\ref{lem:W-properties}, $\E Y_{k,n}^2\to\exp(\sum_{i\ge1}\lambda_i\delta_i^2)=\E(W^{\mathrm{pair}})^2$. For $K>0$, the function $y\mapsto(y-K)_+$ is continuous, and $y\le2(y-K)_+$ for $y>2K$. Therefore,
\begin{align*}
\limsup_{n\to\infty}\E\bigl[Y_{k,n}^2\mathbf 1\{Y_{k,n}^2>2K\}\bigr]
&\le2\limsup_{n\to\infty}\E\bigl(Y_{k,n}^2-K\bigr)_+\\
&=2\lim_{n\to\infty}\Bigl(\E Y_{k,n}^2-\E\bigl(Y_{k,n}^2\wedge K\bigr)\Bigr)
=2\,\E\bigl((W^{\mathrm{pair}})^2-K\bigr)_+ ,
\end{align*}
which tends to $0$ as $K\to\infty$. Thus $(Y_{k,n}^2)_n$ is uniformly integrable, and so is $(Y_{k,n})_n$.

\medskip\noindent
\emph{Step 3: conditioning on simplicity.}
Recall that $\cS=\{X_1=X_2=0\}$, $\Prob(\cS)\to e^{-\lambda_1-\lambda_2}>0$, $G(P)\sim G(n,d)$ conditional on $\cS$, and $Z_k(P)=Z_k(G(P))$ on $\cS$. Let $h(x_1,x_2):=(1-|x_1|-|x_2|)_+$, a bounded continuous function. Since $X_1,X_2$ are integer valued, $h(X_1,X_2)=\mathbf 1_{\cS}$. For any bounded continuous $f$ on $[0,\infty)$ and any bounded continuous $\varphi$ on $\mathbb R^{m-2}$, the continuous mapping theorem and~\eqref{eq:joint-pairing} give
\begin{align*}
\E\bigl[f(Y_{k,n})\varphi(X_3,\dots,X_m)\mathbf 1_{\cS}\bigr]
&=\E\bigl[f(Y_{k,n})\varphi(X_3,\dots,X_m)h(X_1,X_2)\bigr]\\
&\longrightarrow
\E\bigl[f(W^{\mathrm{pair}})\varphi(\xi_3,\dots,\xi_m)\mathbf 1\{\xi_1=\xi_2=0\}\bigr].
\end{align*}
On the event $\{\xi_1=\xi_2=0\}$, the factors of $W^{\mathrm{pair}}$ for $i=1,2$ are deterministic. Since $\lambda_1\delta_1=-\tfrac12$ and $\lambda_2\delta_2=\tfrac14$,
\[
W^{\mathrm{pair}}=e^{-\lambda_1\delta_1-\lambda_2\delta_2}\,W=e^{1/4}\,W\qquad\text{on }\{\xi_1=\xi_2=0\},
\]
where $W$ and $(\xi_3,\dots,\xi_m)$ are independent of $(\xi_1,\xi_2)$. Divide by $\Prob(\cS)\to\Prob(\xi_1=\xi_2=0)=e^{-\lambda_1-\lambda_2}$ to obtain
\begin{equation}\label{eq:cond-distribution}
\mathcal L\bigl(Y_{k,n},X_3,\dots,X_m\bigm|\cS\bigr)\;\Longrightarrow\;\mathcal L\bigl(e^{1/4}W,\xi_3,\dots,\xi_m\bigr).
\end{equation}

The same argument applies to the unbounded functions $f(y)=y$ and $f(y)=y^2$, with $\varphi\equiv1$, by Step~2. Indeed, the continuous mapping theorem gives $Y_{k,n}\mathbf 1_{\cS}\Rightarrow W^{\mathrm{pair}}\mathbf 1\{\xi_1=\xi_2=0\}$ and $Y_{k,n}^2\mathbf 1_{\cS}\Rightarrow(W^{\mathrm{pair}})^2\mathbf 1\{\xi_1=\xi_2=0\}$. These sequences are dominated by the uniformly integrable sequences $Y_{k,n}$ and $Y_{k,n}^2$, respectively, so their expectations converge. Lemma~\ref{lem:W-properties} gives $\E W=1$ and $\E W^2=\exp(\sum_{i\ge3}\lambda_i\delta_i^2)$. Hence
\begin{equation}\label{eq:cond-moments}
c_{k,n}:=\E\bigl[Y_{k,n}\bigm|\cS\bigr]\longrightarrow e^{1/4},
\qquad
\E\bigl[Y_{k,n}^2\bigm|\cS\bigr]\longrightarrow e^{1/2}\exp\Bigl(\sum_{i\ge3}\lambda_i\delta_i^2\Bigr).
\end{equation}

\emph{Proof of \textup{(i)}.}
Since $\E_{G(n,d)}[Z_k(G)]=\E[Z_k(P)\mid\cS]=c_{k,n}\,\E Z_k(P)$, part~(i) follows from~\eqref{eq:cond-moments} and~\eqref{eq:EZk-pairing}.

\emph{Proof of \textup{(ii)}.}
By~\eqref{eq:cond-moments},
\[
\frac{\E_{G(n,d)}[Z_k(G)^2]}{(\E_{G(n,d)}[Z_k(G)])^2}=\frac{\E[Y_{k,n}^2\mid\cS]}{c_{k,n}^2}\longrightarrow\exp\Bigl(\sum_{i\ge3}\lambda_i\delta_i^2\Bigr)
=\exp\Bigl(-\frac{1}{2(d-1)}-\frac{1}{4(d-1)^2}\Bigr)\sqrt{\frac{d-1}{d-2}},
\]
where we used $\lambda_1\delta_1^2=\frac{1}{2(d-1)}$, $\lambda_2\delta_2^2=\frac{1}{4(d-1)^2}$, and $\frac{1}{2(d-1)}+\frac{1}{4(d-1)^2}=\frac{2d-1}{4(d-1)^2}$.

\emph{Proof of \textup{(iii)}.}
Conditional on $\cS$, $Z_k(G)/\E_{G(n,d)}[Z_k(G)]=Y_{k,n}/c_{k,n}$, where $c_{k,n}\to e^{1/4}$. Thus~\eqref{eq:cond-distribution} and Slutsky's theorem give the joint convergence $(Z_k(G)/\E_{G(n,d)}[Z_k(G)],X_3(G),\dots,X_m(G))\Rightarrow(W,\xi_3,\dots,\xi_m)$ for $G\sim G(n,d)$. Lemma~\ref{lem:W-properties} gives the stated properties of $W$.

\medskip\noindent
\emph{Step 4: a common approximation.}
Let $W_n^{(m)}:=\prod_{i\le m}(1+\delta_i)^{X_i}e^{-\lambda_i\delta_i}$. This is a function of the pairing $P$ alone and does not depend on $k$. The map $(y,x_1,\dots,x_m)\mapsto y-\prod_{i\le m}(1+\delta_i)^{x_i}e^{-\lambda_i\delta_i}$ is continuous. Therefore~\eqref{eq:joint-pairing} and the continuous mapping theorem give $Y_{k,n}-W_n^{(m)}\Rightarrow W^{\mathrm{pair}}-W^{(m)}(\xi)$. By the portmanteau theorem, for every $\varepsilon>0$,
\[
\limsup_{n\to\infty}\Prob\bigl(|Y_{k,n}-W_n^{(m)}|\ge\varepsilon\bigr)\le\Prob\bigl(|W^{\mathrm{pair}}-W^{(m)}(\xi)|\ge\varepsilon\bigr)\xrightarrow[m\to\infty]{}0,
\]
because $W^{(m)}(\xi)\to W^{\mathrm{pair}}$ almost surely. Thus
\begin{equation}\label{eq:common-approx}
\lim_{m\to\infty}\limsup_{n\to\infty}\Prob\bigl(|Y_{k,n}-W_n^{(m)}|\ge\varepsilon\bigr)=0\qquad\text{for every }\varepsilon>0\text{ and every }1\le k\le d-1 .
\end{equation}

\medskip\noindent
\emph{Step 5: proof of \textup{(iv)}.}
Let $1\le k,\ell\le d-1$, and let $n$ be such that $dn$, $kn$, and $\ell n$ are even. Then $Y_{k,n}$ and $Y_{\ell,n}$ are defined on the same space $\mathcal P_{n,d}$. For $\varepsilon>0$ and every $m$,
\[
\Prob\bigl(|Y_{k,n}-Y_{\ell,n}|\ge2\varepsilon\bigr)\le\Prob\bigl(|Y_{k,n}-W_n^{(m)}|\ge\varepsilon\bigr)+\Prob\bigl(|Y_{\ell,n}-W_n^{(m)}|\ge\varepsilon\bigr);
\]
taking $\limsup_{n\to\infty}$ and then $m\to\infty$ in~\eqref{eq:common-approx} gives $Y_{k,n}-Y_{\ell,n}\to0$ in probability in $\mathcal P_{n,d}$. Since $\Prob(A\mid\cS)\le\Prob(A)/\Prob(\cS)$ and $\Prob(\cS)$ is bounded away from zero, also $Y_{k,n}-Y_{\ell,n}\to0$ in probability conditional on $\cS$, and hence for $G\sim G(n,d)$. Now, for $G\sim G(n,d)$,
\[
\frac{Z_k(G)}{\E_{G(n,d)}[Z_k(G)]}-\frac{Z_\ell(G)}{\E_{G(n,d)}[Z_\ell(G)]}
=\frac{Y_{k,n}}{c_{k,n}}-\frac{Y_{\ell,n}}{c_{\ell,n}}
=\frac{Y_{k,n}-Y_{\ell,n}}{c_{k,n}}+Y_{\ell,n}\Bigl(\frac1{c_{k,n}}-\frac1{c_{\ell,n}}\Bigr).
\]
The first term tends to $0$ in probability because $c_{k,n}\to e^{1/4}$. The second does likewise because $Y_{\ell,n}$ is tight, by part~(iii), while $1/c_{k,n}-1/c_{\ell,n}\to0$ by~\eqref{eq:cond-moments}. This proves the first assertion in~\eqref{eq:linked}.

For the second assertion, let $D_n$ be the difference just shown to be $o_p(1)$, and let $V_n:=Z_\ell(G)/\E_{G(n,d)}[Z_\ell(G)]$. Then $V_n\Rightarrow W$, where $W>0$ almost surely by part~(iii). Given $\eta>0$, choose $\varepsilon>0$ such that $\Prob(W\le\varepsilon)\le\eta$ and $\varepsilon$ is a continuity point of the distribution of $W$. Then $\limsup_n\Prob(V_n\le\varepsilon)\le\eta$. On $\{V_n>\varepsilon\}$, we have $Z_\ell(G)>0$ and
\[
\Bigl|\frac{Z_k(G)/\E_{G(n,d)}[Z_k(G)]}{Z_\ell(G)/\E_{G(n,d)}[Z_\ell(G)]}-1\Bigr|=\frac{|D_n|}{V_n}\le\frac{|D_n|}{\varepsilon}\xrightarrow{\;p\;}0 .
\]
Since $\eta>0$ is arbitrary, the ratio converges to $1$ in probability. This completes the proof of Theorem~\ref{thm:main}. \qed

\begin{remark}
The argument gives a stronger form of~\eqref{eq:linked}. By~\eqref{eq:common-approx}, for every $k$, the normalised count $Z_k(G)/\E_{G(n,d)}[Z_k(G)]$ is approximated in probability, as $n\to\infty$ followed by $m\to\infty$, by the explicit function $e^{-1/4}W_n^{(m)}=\prod_{i=3}^m(1+\delta_i)^{X_i(G)}e^{-\lambda_i\delta_i}$ of the short cycles of $G$. On simple graphs, $X_1=X_2=0$, and the factors for $i=1,2$ in $W_n^{(m)}$ reduce to the constant $e^{1/4}$. The approximation does not depend on $k$ because the transfer matrix $T$ in Lemma~\ref{lem:short-cycles} has the same spectrum $\{d-1,-1\}$ for every $k$, although the matrix itself depends on $k$.
\end{remark}

{\small
\let\oldthebibliography\thebibliography
\renewcommand{\thebibliography}[1]{\oldthebibliography{#1}\setlength{\itemsep}{0pt}\setlength{\parskip}{1pt}}

}

\end{document}